\pdfoutput=1
\documentclass{article}
\usepackage[margin=25mm]{geometry}
\usepackage{amsmath}
\usepackage{amsfonts}
\usepackage{amssymb}
\usepackage{graphicx}
\usepackage[utf8]{inputenc}
\usepackage[T2A]{fontenc}
\usepackage[ukrainian,russian,english]{babel}
\usepackage[unicode, pdftex]{hyperref}
\usepackage{xcolor}
\usepackage[all]{xy}
\newtheorem{lemma}{Lemma}
\newtheorem{theorem}{Theorem}
\newtheorem{definition}{Definition}
\newtheorem{proof}{Proof}

\title{Structure of optimal gradient flows bifurcations on closed surfaces}

\author{Illia Ovtsynov, Alexandr Prishlyak} 

\begin{document}

\maketitle
\begin{abstract}
  We consider structure of typical gradient flows bifurcations on closed  surfaces with minimal number of singular points. There are two type of such bifurcations: saddle-node (SN) and saddle connections (SC). The structure of a bifurcation is determinated by codimension one flow in the moment of bifurcation.  We use  the chord diagrams to specify the flows up to topological  trajectory equivivalence.  A chord diagram with a marked arc is complete topological invariant of a SN-bifurcations and a chord diagram with T-insert -- of SC-bifurcations.  We list all such diagrams for flows on norientable surfaces of genus at most 2 and nonorientable surfaces of genus at most 3. For each of diagram we found inverse one that correspond the inverse flow.
\end{abstract}

\section{Introduction}
Let's suppose that $M$ is a closed 2-manifold. Two flows $\{g^t_i: M \rightarrow M\}$, $i = 1, 2$, are said to be \textit{topologically or trajectory equivalent}, if there is a homeomorphism $h: M \rightarrow M$
which maps trajectories of  $g^t_1$ onto trajectories of $g^t_2$.

Two 1-parametrical families of flows $g^t_s, G^t_S, \ s, S\in [0,1]$,   are\textit{ topologically equivalent} if there is such a parameter replacement $s\to S$ that creates a family of topologically equivalent flows (vector fields). \textit{Typical family} of flows (vector fields) is such a family which remains topologically equivalent to the initial one after small change of the parameters defining it (as family of vector fields).

Codimension zero flows on surfaces are Morse--Smale flows.
\textit{Codimension one flows} are flows which appear in typical 1-parametrical families of flows and which are not codimension zero flows.

One can build complete topological invariants for topological classification of flows. They are often graphs having special properties and being endowed with additional information. Wherein, two flows are topologically equivalent if and only if there is an isomorphism between the graphs which preserves this additional information. Effectiveness of invariants is defined with an opportunity to find all the structures of dynamical systems with a given number of singular points.


For flows with a finite codimension (they have a finite number of singular points, closed trajectories and separatrices), we have a separatrix diagram as a complete topological invariant \cite{Fleitas1975, Leontovich1955,Oshemkov1998, Peixoto1973}.

We can simplify this invariant for Morse flows (Morse--Smale flows without closed trajectories) which are topologically equivalent to gradient flows of Morse function. It is a graph embedded into a surface. Its vertices are sources and its edges are stable manifolds of saddles.

Rotation system \cite{GT87} is often used for defining graph embedding into a surface. This is, in fact, a Peixoto invariant \cite{Peixoto1973} for Morse flows. In the case of one-source flow it is convenient to use chord diagrams for its defining \cite{Kybalko2018}.


Also, complete topological invariants for flows with gradient dynamics were constructed in dimension 2 \cite{prishlyak2020three, Prishlyak2022} and dimension 3 \cite{Prishlyak2002beh2, Hatamian2020, prishlyak2003sum, prish1998, prish2001, Prishlyak2023}, as well as for collective dynamics in dimension 2  \cite{Prishlyak2020, Prishlyak2021}.


Besides, to illustrate some bifurcations of singular points,
transitions via heteroclinic separatrices and so on, one can characterize and list all generic
non-Morse gradient flows \cite{Kibkalo2022}.

The the structures of bifurcations are similar to the structure of function deformations~\cite{Bilun2002, hladysh2019simple, Hladysh2019, Khohliyk2020, Kravchenko_2019, Maks2006, Maks2020, Prishlyak2000, prishlyak2001conjugacy,Prishlyak2023a, Reeb1946, prish2002tm, prish2015}. There the codimension zero functions are Morse functions of general position. Their structure is defined with Reeb graphs, their deformations are defined with deformations of Reeb graphs. Codimension one functions have a cubical speciality or two critical points in the same level. However, if the bifurcation ``saddle--node{}'' can be described as a gradient of a deformation with a cubical speciality, then the saddle connection is not associated with two critical points in the same function level.

The purpose of this paper is constructing complete invariant for codimension one flow on closed surfaces which resembles a chord diagram for Morse flows. This invariant has a marked point in the diagram, which allows us to define clearly a number code of the flow.

In the section~\ref{prel} we consider Morse--Smale flows on 2-manifolds, what optimal Morse flows are and what complete topological invariants are used to describe these flows. It helps reader to understand better the description of optimal gradient codimension one flows in this article. Also we explain what the reverse diagram is when we consider SN-flows and SC-flows and how to build reverse codes.

In the section~\ref{main} we research gradient codimension one flows on closed 2-manifolds which can be divided into two types --- SN-flows (saddle--node flows) and SC-flows (saddle connection flows). We speak about the appropriate topological invariants for such flows. 

In the section~\ref{prove} we place all the proofs of relevant theorems. 

In the section~\ref{examples} we show all the examples of SN-flows and SC-flows on both oriented and non-oriented surfaces.

\section{Preliminaries}\label{prel}

Let's suppose that $M$ is a smooth surface (2-manifold) and $X$ is a smooth vector field on this surface which generates the flow $g^t$.




\begin{definition}
A \textit{stable manifold} (denoted by $W^s(p)$ or $S(p)$) of a 
singular point $p$ is the following set
\begin{equation*}
  S(p) = \{x \in M: \lim\limits_{t \rightarrow \infty} g^t(x) = p\}.
\end{equation*}
Besides, the set
\begin{equation*}
  U(p) = \{x \in M: \lim\limits_{t \rightarrow -\infty} g^t(x) = p\}
\end{equation*}
is called an \textit{unstable manifold} (denoted by $W^u(p)$ or $U(p)$) of $p$.
\end{definition}

\begin{definition}
A point $p \in M$ is called \textit{non-wandering} for the vector field $X$ if for any $T \in \mathbb{R}$ and any neighbourhood $U$ of $x$ there exists $t > T$ such that $g^t(U) \cap U \neq \varnothing$. Otherwise the point is said to be \textit{wandering}.
\end{definition}

The set of non-wandering points of the vector field $X$ is denoted by $\Sigma(X)$.

\begin{definition}
A vector field $X$ (flow $g^t$) on a compact manifold $M$ is a \textit{Morse vector field (flow)}, or \textit{gradient Morse--Smale vector field (flow)} if:
\begin{enumerate}
  \item[\rm 1)] $X$ has a finite number of critical points and all of them are hyperbolic;
  \item[\rm 2)] the set of non-wandering points $\Sigma(X)$ coincides with the union of all 
singular points of $X$;
  \item[\rm 3)] if $S(p)$ and $U(q)$ are stable and unstable manifolds of singular points of $p$ and $q$, and $x \in S(p) \cap U(q)$, then $T_x S(p) + T_x U(q) = T_x M$, i.e. the manifolds $S(p)$ and $U(q)$ intersect at $x$ transversally.
\end{enumerate}
\end{definition}

Item 2) in the definition of Morse vector field can be replaced by this statement: \textit{each trajectory of the flow has its alpha- and omega-limit sets in singular points}.

Topological equivalence of optimal flows, defined on 2-manifolds, often can be researched via using the notion of chord diagrams.

\begin{definition}
Let's suppose that $n \geq 0$. A configuration consisting of a circle, $2n$ distinct points on it, and $n$ chords, which specify a partition of these points into pairs, will be called \textit{a chord diagram with n chords}, or \textit{n-diagram}.
\end{definition}

\begin{definition}
Two chord diagrams are said to be \textit{equivalent} if there is a
homeomorphism of circle to circle that preserves partition of the points into pairs (ends of chords). If the homeomorphism preserves the orientation of the circle, then we call such diagrams \textit{isomorphic}.
\end{definition}

\begin{definition}
A Morse flow on a surface is called \textit{optimal} if it has the least number of singularities among all Morse flows on this surface.
\end{definition}



Suppose that $X$ is an optimal Morse field on a connected surface. Let's take such a small neighbourhood $U$ for the source $x_0$ that its boundary $\partial U(x_0)$ is a smooth curve and it is transversal to each flow trajectory.

Let's take a stable manifold for each saddle that consists of all the trajectories converging to the saddle. These trajectories start at $x_0$, so, they intersect $\partial U(x_0)$ at one point. Afterwards we get a couple of points on the boundary which is homeomorphic to a circle $S^1 = \{(x, y) | x^2 + y^2 = 1\}$.

It means that the couple of points maps into a couple of points on the circle $S^1$. Connecting the points of the same couple with the chord we
will obtain chords at the circle, which do not have common ends, for the
stable manifolds of the saddle points do not intersect with each other. Thus
we obtain a chord diagram. If the union of $U$ with a regular neighbourhood of
the stable manifold is an oriented surface, then
we leave chord unmarked. Otherwise we get it marked (we just put a symbol X on this chord).

We are drawing reader's attention to the fact that such diagrams can be drawn only for those Morse fields which have only one source.

\begin{definition}
A \textit{marked chord diagram} is a chord diagram with marking of its chord. If the homeomorphism preserves the orientation of the circle, then we call such diagrams \textit{isomorphic}.
\end{definition}

\begin{definition}
Marked chord diagrams are said to be \textit{equivalent} if there is a
homeomorphism of circle to circle that preserves the whole structure of the circle, including chords and their marking.
\end{definition}

The examples of chord diagrams can be seen in the work~\cite{Kybalko2018}, at the figures 8.1--8.6.

A marked chord diagram is a complete topological
invariant of the optimal Morse vector field. Since now we call it a \textit{chord diagram of a
vector field}.

Let's choose a point on a
certain arc that is not the end of a chord and the counterclockwise direction of circumvention
of a circle. We move from the chosen point
in this direction to the nearest end of the chord. At every such point the
movement on the circle changes on the chord movement. After passing the
chord we continue the circular motion in counterclockwise direction if the
chord is unmarked and change the direction along to clockwise if the chord
is marked. We move to the nearest point that is the end of the chord. In it
again change the movement on the movement of the chord. After marked
chords we change the direction and after unmarked one we preserve it. We will continue this process until we return to the starting
point. If you want to research these structures in more details, see the article~\cite{Kybalko2018}, particularly the figure 6.1 which illustrates the diagram we have just spoken about.

We obtained a cycle which we call a \textit{chord cycle}. As soon as chord cycle covers all chords (two times) and arcs we have \textit{a one-cycle diagram}. Chord diagram of optimal Morse flow is one-cycle~\cite{Kybalko2018}.




\section{Main results} \label{main}
All the proofs of the theorems given in this section can be read in the chapter~\ref{prove}.

One of the conditions in the definition of Morse--Smale flow is not true for typical one-parameter bifurcations. For gradient flows we have next possible violations: 1) singular point is non-hyperbolic --- saddle--node singularity (SN) and 2) existence of saddle connection (SC). Therefore codimension one gradient flows are divided into two types --- SN and SC.

\subsection*{Saddle--node}
\begin{figure}[ht]
\center{\includegraphics[width=0.8\linewidth]{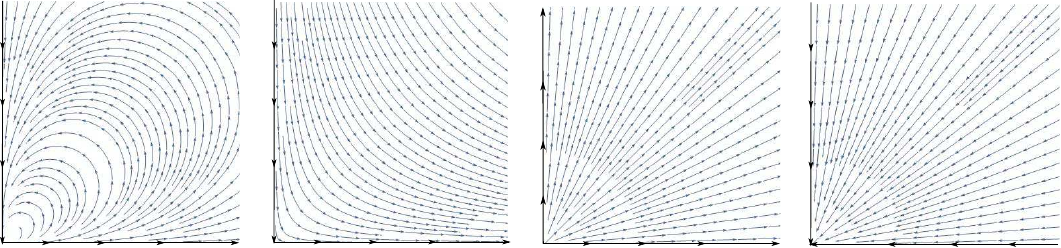}
}
\put (-255,-15) {1}
\put (-180,-15) {2}
\put (-105,-15) {3}
\put (-30,-15) {4}
\caption{Different kinds of angles in the neighbourhood of fixed points}
\label{angles}
\end{figure}
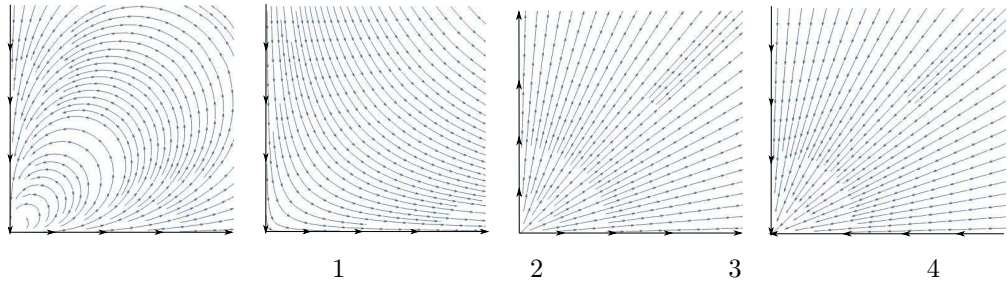

The traectories going out of and into fixed points divide their neighbourhoods into different angles. As a conclusion, we can have different types of these angles depending on how the trajectories in the neighbourhood look. In the figure~\ref{angles} we have an elliptic angle in the picture 1, a hyperbolic angle in the picture 2, a parabolic sink in the picture 3, a parabolic source in the picture 4.

A flow is called \textit{SN-flow} if 1) it has one singular saddle--node type point, but the other singular points are hyperbolic, 2) all the other conditions of Morse flow are true, 3) the separatrices which intersect boundary of saddle--node hyperbolic angle have their alpha- or omega-limit sets in the source or in the sink.

If the stable manifold of a singular saddle--node type consists of only one trajectory, then this point has the type ``saddle--source{}''. If the unstable manifold of a singular saddle--node type consists of only one trajectory, then this point has the type ``saddle--sink{}''.

\begin{figure}[ht]
\center{\includegraphics[width=0.8\linewidth]{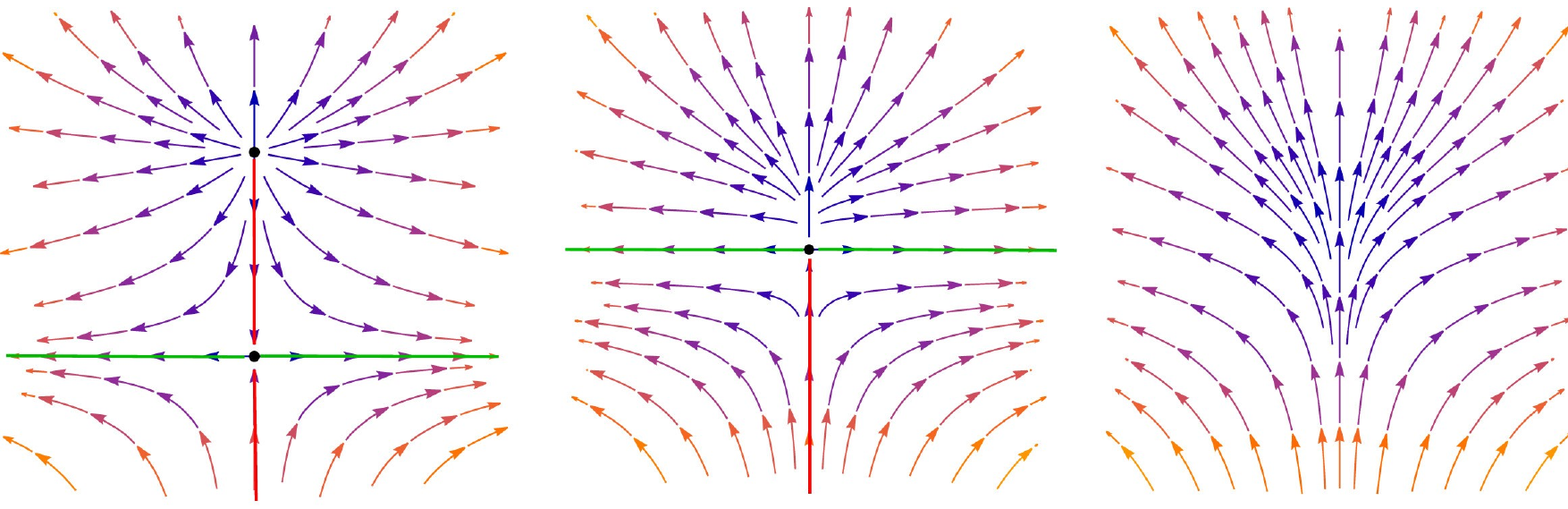}
}
\caption{Saddle--node bifurcation}
\label{exsn}
\end{figure}

We can observe the bifurcation saddle--node when the node is a source in the figure~\ref{exsn}, at the picture 2. It can be described with the equation $V(x, y, a) = \{x, y^2 + a\}$, where $a$ is a parameter. When $a < 0$, we will get the flow before the bifurcation (left picture), when $a > 0$, we will get the flow after the bifurcation (right picture), and when $a = 0$, we will get the flow at the moment of the bifurcation (centre picture).

An SN-flow on a surface M is called \textit{optimal} if there is no other SM-flow with a fewer number of singular points on M.


\begin{lemma}\label{lem1}
SN-flow is optimal if and only if it has one source and one sink.
\end{lemma}

Let's build a C-diagram for an optimal SN-flow which is a chord diagram with a marked chord C. We will consider that a saddle--node singularity is created with a connection between the source and the saddle (the case sink--saddle can be got with reversing the motion along all the trajectories). Let's analyze a separatrix $\gamma$, which is moving into the saddle--node SN. According to the definition of SN-flow, its $\alpha$-limit set is the source $S$. Let's suppose that $L$ is a boundary of a regular neighbourhood $U$ of the union $S \cup \gamma \cup SN$ such that flow in the points in $L$ is directed outward of $U$. Then we will build the chord diagram by analogy with polar Morse flows: we will map $L$ homeomorphically into a unit circle, the images of the intersections $L$ with other separatrices will be connected by chords if they belong to the same stable manifold of the saddle point. The intersection of the unstable manifold of saddle--node with $L$ is mapped into an arc of the unit circle. Hence we have a chord diagram with a marked arc (C-diagram). Let's colour the chords in the diagram. All the chords are coloured with the first colour for an oriented surface. A chord is coloured with the second colour for an non-oriented surface if there is a neighbourhood of the appropriate stable manifold which constructs a M$\ddot{\textrm{o}}$bius strip with the neighbourhood of the source. Otherwise the chord is coloured with the first colour.

The example of C-diagrams is given, for instance, in the figure~\ref{bif1}.

The components of the boundary of the union $U$ with regular neighbourhoods of stable 1-manifolds represent the cycles on the chord diagram. A chord diagram with one cycle is called \textit{one--cycled}. A chord diagram of an optimal flow is one--cycled what implies from the lemma.


\begin{theorem}\label{t1}
  Two optimal SN-flows are topologically equivalent if and only if they have the same type of singularity ``saddle--node{}'' and there is an isomorphism between their chord diagrams which keeps marked arcs.
\end{theorem}


\begin{theorem}\label{t2}
Let's assume that there is a marked arc of a unit circle in a one--cycled chord diagram, ends of which don't coincide with the ends of the chords. Then there is an optimal SN-flow with this chord diagram.
\end{theorem}


Let's fix an orientation on the unit circle of the C-diagram, having defined the counterclockwise moving direction. Then this direction is defining an orientation of the arc C and can be defined after ordering the ends of the arc C. We begin the motion from the first end of C along the unit circle along the orientation for constructing the code of the C-diagram. We also denote the other end of C by 0. We enumerate the chords by the order of passing the first one of their ends which are met during the motion along the circle. Let's write the sequence of numbers from the ends of the chords and the second end of the arc C. So, 0 is present in this sequence for once, whereas the other numbers are present each for twice (both ends of the chords). If a chord is coloured with the second colour (non-oriented) then we put a dash after the first of its numbers ('). We call a diagram code \textit{symmetric} to other ones if there is an axial symmetry for this diagram.

For instance, if C is an upper half-circle in the diagram 4, then the code is following: 1'210323. This code is symmetric.

\begin{theorem}\label{t3}
 Two optimal SN-flows are topologically equivalent if and only if they have either the same or isomorphic codes.
\end{theorem}

\subsection*{Saddle connection}

A flow is called an \textit{SC-flow} if 1) it has a finite number of singular points and they are hyperbolic, 2) each trajectory begins in a singular points and ends in a singular point, 3) this flow has a unique separatrix connection.

\begin{figure}[ht]
\center{\includegraphics[width=0.8\linewidth]{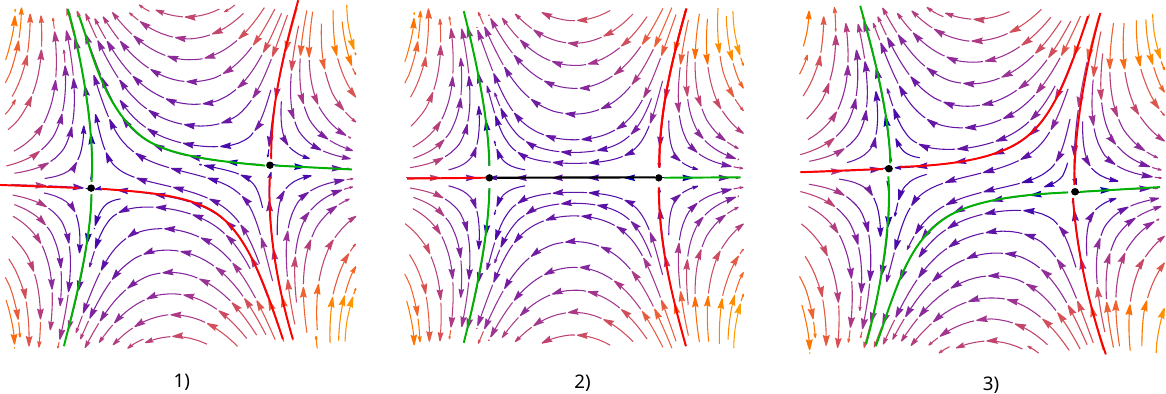}
}
\caption{Saddle connection bifurcation}
\label{exsc}
\end{figure}

The examples of saddle connections are given in the figure~\ref{exsc}.

SC-flow on surface M is called \textit{optimal} if there is no other flow with a fewer number of singular points on M.

\begin{lemma}\label{L2}
An SC-flow on a surface of nonpositive Euler characteristics is optimal if and only if it has a unique source and a unique sink.
\end{lemma}

Let's call \textit{T-graph} a connected graph with four vertices where three ones of them have degree 1 but another one has degree 3. One of the three first vertices is marked and called \textit{lower} vertice whereas two others are called \textit{side} vertices. Correspondent adjacent edges to them are called lower and side ones.

Let's build a T-diagram of an optimal SC-flow which is a chord diagram with inscribed T-graph in it so that all ends of the graph are on the unit circle and are not ends of any chords. As for saddle points there is no a connection between any of them, we construct a chord diagram for them as for a polar Morse flow. Let's aplha-limit of a saddle connection is $x$ and omega-limit is $y$. We build for $x$ a standard chord $L$, which represents the stable manifold. Let's denote a middle of this chord by $z$. A stable manifold of $y$ represents a point in a unit circle. We join this point with $z$ with a segment $N$. $T = L \cup N$ be called a \textit{T-graph}, $N$ is lower edge, halves of $L$ ($L_1$,$L_2$) are side edges. A chord diagram with a T-graph is called a \textit{T-diagram}. Besides, all the chords, $N$, $L_1$ and $L_2$ are coloured with two colours.


Let's suppose that a simple regular curve $\gamma$ is transversally intersected with oriented curves $\alpha$ and $\beta$ and $x=\alpha \cap \gamma \ne \beta \cap \gamma =y$. We  say that orientations of the curves $\alpha$ and $\beta$ are coherent if there is a simply connected neighbourhood $U$ of the arc $\gamma$ between $x$ and $y$ such that $\alpha \subset \partial U$, $\beta \subset \partial U$ and the orientations of $\alpha$ and $\beta$ represent the same orientation $\partial U$. A chord in the chord diagram is coloured with the first colour if an oriented path around the boundary of the neighbourhood of the source in neighbourhoods of the points of intersection with the appropriate separatrix represents coherent orientations. Otherwise the chord is coloured with the second colour.


Let's suppose that $x$ is alpha-limit and $y$ is omega-limit of the saddle-connection $\gamma$, $U$ is a regular neighbourhood of the source, $V$ is a regular neighbourhood of $x \cup \gamma \cup y$. Let's denote $a=W^s(y) \cap \partial U$,  $b \cup c=W^s(x) \cap \partial U$, $a_1=W^s(y) \cap \partial V$,  $b_1 \cup c_1=W^s(x) \cap \partial V$, $b$ and $b_1$ belong to the same separatrix. If the orientation of $\partial U$ is defined with the ordered triple $\{ a,b,c\}$, then the orientation of $\partial V$  be defined with the ordered triple $\{ a_1,c_1, b_1\}$. Then we have two oriented transversal intersections with $\partial U$ and $\partial V$ for each among three separatrices which go through the points $ a,b,c$. If their orientations are coherent relatively to the separatrix, then we colour  the appropriate edge with the first colour (leave unmarked), otherwise --- with the second colour (mark it with X on  the diagram).


Examples of such diagrams can be seen in the figure~\ref{snn}.

\begin{theorem}\label{t4}
  Two optimal SC-flows are topologically equivalent if and only if their T-diagrams are isomorphic.
\end{theorem}


\begin{theorem}\label{t5}
  If we have a T-diagram, then there is an optimal SC-flow with this diagram.
\end{theorem}


\begin{theorem}\label{t6}
  A T-diagram for optimal SC-flows on oriented surfaces and only on them has following properties: 1) all the chords are coloured with the first colour, 2) all the edges in T have the same colour.
\end{theorem}


We enumerate the lower end of T as 0 and the other ones as 1 for defining a code of a T-diagram. The chords are enumerated depending on the order of their first end we meet during the counterclockwise motion along the unit circle beginning from 0. The code is a sequence of numbers of the ends of the chords and the edges which we meet during this motion. We put the dash ' after the numbers of non-oriented edges and the first number of a non-oriented chord (coloured with the second colour). We can get a symmetric code if we do the same doing clockwise motion.


For example, there is a code 01'212 for the diagram 6 in the fig. \ref{T3a}, the symmetric code is 02121'.


\begin{theorem}  \label{t7}
Two optimal SC-flows are topologically equivalent if and only if they have either the same or isomorphic codes.
\end{theorem}


\subsection*{Reverse diagram and reverse code}

A chord diagram of an optimal SN-flow be reverse for itself, only the type of the singular point saddle--node be changed.

Let's analyze the construction of the reverse T-diagram. We orientate the cycle which is a boundary of the union of the half-sphere with a neighbourhood of T and unintersected strips considered as neighbourhoods of the chords. Considering these strips as ones which are homeomorphic to $[0, 1] \times [0, 1]$, we can define transversal to them curves which are homeomorphic to $t \times [0, 1]$. These curves intersect the cycle twice. Besides, a strip $[0, 1] \times [0, 1]$ is englued to the half-sphere by the embedding $\{-1, 1\} \times [-1, 1]$ in a unit circle, whereas the chord corresponds to the set $[-1, 1] \times 0$. Let's map the cycle into the unit sphere and join the appropriate points with chords. A colour of a chord can be defined by coherence of motion directions along the cycle with relation to the transversal intersection. Continuing the lower edge of T until the intersection with the cycle create a basis of T in the reverse diagram, whereas two other ends of reverse T we get as intersection of the transversal and the cycle before the middle of the lower edge of initial T (1a and 1b).

\begin{figure}[ht]
\center{\includegraphics[height=3.3cm]{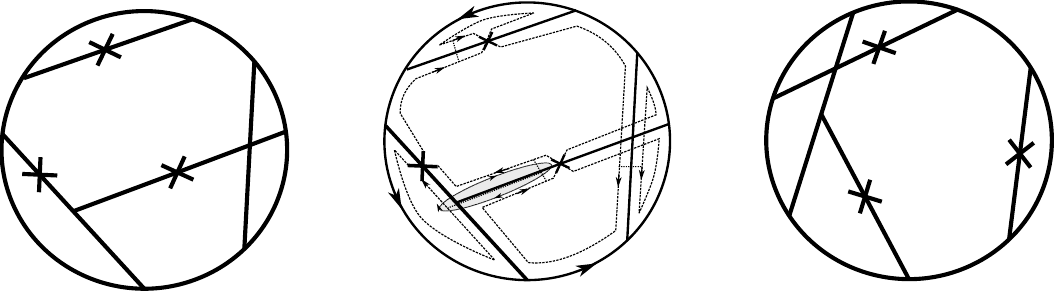}}
\put (-204,20) {\small{0}}
\put (-175,45) {\small{1a}}
\put (-165,27) {\small{1b}}
\put (-193,65) {\small{2}}
\put (-200,87) {\small{2}}
\put (-147,34) {\small{3}}
\put (-132,36) {\small{3}}
\put (-5,75) {\small{2}}
\put (-15,5) {\small{2}}
\put (-28,93) {\small{3}}
\put (-69,93) {\small{1b}}
\put (-100,55) {\small{3}}
\put (-91,13) {\small{1a}}
\put (-47,-7) {\small{0}}
\caption{T-diagram (left), dotted cycle on the diagram (center), reverse diagram (right) }
\label{rev}
\end{figure}
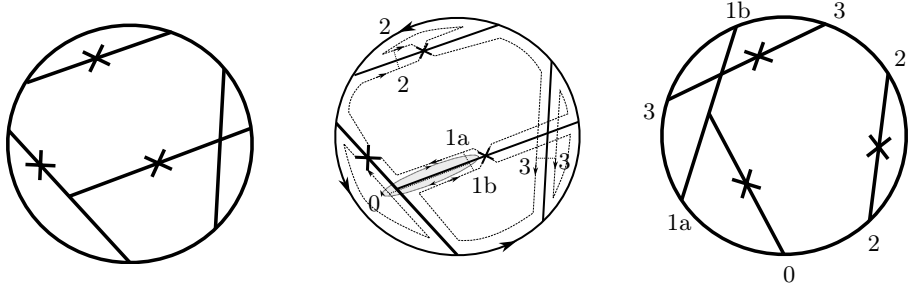

Orientation of the boundary of the neighbourhood of new T is defined by analogy with initial T: we put a bijection between an ordered triple of points in cycle (the basis of T and two other ends of T) and a triple of points in the neighbourhood of reverse T and reverse the order in the triple. If the orientation of the cycle and boundaries of the neighbourhood are coherent in the appropriate edges of reverse T then these edges are coloured with the first colour, otherwise --- with the second colour. We call half-edge of the lower edge of T from its center to the center of lower edge central \textit{half-edge}. Let's assume that orientation in the unit circle of the chord diagram is defined with the counterclockwise motion. Then the orientation on the boundary of the neighbourhood of the central half-edge which is coherent with the orientation on the unit circle is defined with the clockwise motion. If this orientation is coherent (we mean opposite) with the motion direction along the cycle in points 0, 1a and 1b, then the appropriate edge is coloured with the first colour, otherwise --- with the second one.

\begin{theorem}\label{t8}
  The diagram which is reverse to the T-diagram of flow T, is a T-diagram of reverse flow $-X$.
\end{theorem}

\section{Proofs of theorems} \label{prove}

\begin{proof} (lemma \ref{lem1}) Let's $R$ is a saddle-node of source type.
According to the definition of SN-flow, if $\gamma$ is a separatrix and $\omega(\gamma)=R$, then $\alpha(\gamma)$ is a source. If $\gamma$ is a separatrix and $\alpha(\gamma)=R$, then $\omega(\gamma)$ is a sink.  So we have a source and a sink in the SN-flow. On the other hand, there is always an SN-flow with one source and one sink on a surface which we call a \textit{polar flow}. We can get this flow, making a bifurcation saddle--node in any regular point of an optimal Morse flow. If a flow has more than one source or more than one sink, then according to the Poincar$\acute{\textrm{e}}$--Hopf theorem a number of saddles is larger than the same number for a polar flow too. Hence a polar flow is an optimal one.
\end{proof}

\begin{proof} (theorem \ref{t1}) Necessity. A topological equivalence of flow maps the stable manifolds of one flow into the stable ones of the other one. It defines a correspondence between the chords and the ends of marked arcs, hence, there is an isomorphism between the C-diagrams. 

Sufficiency. If there is an isomorphism between C-diagrams, then it defines a correspondence between the appropriate separatrices. Each of the flows has the only a separatrix ($\gamma$) which alpha-limit is the source and omega limit is the saddle--node. So we can prolong this correspondence between the separatrices on $\gamma$. A cyclical order of the separatrices in the source is defined by the order of the ends of the chords which do not belong to the C-arc. A cyclical order of the separatrices which go out of the saddle--node is defined with the order of the ends of the chords, which belong to the C-arc, and the ends of the C-arc. A cyclical order of the separatrices in the saddles and in the sink is defined with the C-diagram too. So, the separatrix diagrams of the flows are isomorphic and the flows are topologically equivalent according to the Peixoto classification \cite{Peixoto1973, Oshemkov1998}.
\end{proof}

\begin{proof} (theorem \ref{t2}) Let's build a graph with two vertices which correspond to the sink and the saddle--node. Edges of the graph correspond stable manifolds of saddle points and the saddle--node. Like in the proof of the previous theorem, the C-diagram defines cyclical orders of half-edges in vertices. Then according to the topological graph theory there is a cell embedding (with one 2-cell which corresponds to the cycle of the C-diagram) of the built graph into a closed surface. This embedding defines a separatrix diagram with an appropriate flow for it.
\end{proof}

\begin{proof} (theorem \ref{t3})
If we fix an orientation in the unit circle of a C-diagram, then according to the procedure of constructing a circle, the code is defined clearly. For the other orientation we have a symmetric code. So codes are either identical or symmetric for topologically equivalent flows. On the other hand, we can renew the C-diagram clearly using the code, if we have a regular 2n+2-gon with its vertices in the unit circle (n is a quantity of mutually distinct numbers in the code apart from 0) and want to put numbers on its vertices according to the code (doing a clockwise motion along the circle), whereas the second 0 be put on the last vertice. The pairs of points with identical numbers (apart from 0) be joined by chords. The arc which has been moved from the second zero to the first zero in a counterclockwise way be marked as a C-arc.
\end{proof}

\begin{proof} (lemma \ref{L2})
The proof is coonstructed with the analogy to the proof of the lemma \ref{lem1}.
\end{proof}

\begin{proof} (theorem \ref{t4})
Necessity. Topological equivalence of flows defines a correspondence between their trajectories, hence, a correspondence between their intersections with boundary of regular neighbourhood of source. So, it defines a homeomorphism of unit circles of the chord diagrams which preserves the ends of the chords and the edges, so it is an isomorphism of the T-diagrams.

Sufficiency. Every T-flow has a unique separatrix connection, so topological equivalence maps a separatrix connection into another separatrix connection. The mapping of other stable 1-manifolds (separatrices) is defined with the correspondence between the chords and the edges of the T-diagrams. An order of the ends of the chords and the edges in the unit circle of the T-diagram defines a cyclical order of the separatrices in the source, so separatrix diagrams of the flows are isomorphic, whereas the flows themselves are topologically equivalent.
\end{proof}

\begin{proof} (theorem \ref{t5})
Using the T-diagram we construct a separatrix diagram of the flow. For this we englue a 0-handle which is a hemisphere to a unit circle (boundary of the hemisphere), each chord corresponds to a strip (1-handle), englued to a hemisphere. It is twisted if a colour of the chord is the second one and it isn't if a colour of the chord is the first one. 1-handles do not intersect each other. We take 2-disc of smaller radius with a center in a triple point T in a plane inside the circle. Let's orientate its boundary using clockwise motion. Three edges of the T-diagram correspond to three strips (1-handles), which join this 2-disc with a hemisphere. Their twisting is defined with a colour of the edges, like in the case with the chords. These three 1-handles together with 2-disc construct T-neighbourhood. Each 1-handle has an axis (middle edge). We join its ends with the shortest arc with a lower point in the hemisphere if it lies in the unit circle, or with the segment with its center in 2-disc otherwise. We fix points in each loop which corresponds a chord, lower edge of T and the center of T --- these points correspond saddles. The source corresponds to the lower point of the sphere. Prolonged axes of the 1-handles correspond to stable separatrices. Let's englue 2-disc (2-handle) to the boundary of the created surface. Coaxes intersect this boundary in the points which correspond unstable separatrices. We join them with the center of the 2-handle. Here we got unstable separatrices. So we have a closed 2-manifold and a separatrix diagram on it which defines a flow.
\end{proof}

\begin{proof} (theorem \ref{t6})
If there were a chord of another colour in the diagram, then the appropriate 1-handle with the 1-handle would construct a M$\ddot{\textrm{o}}$bius strip, but that is impossible on an oriented surface. As well, if all the T-edges are of the first colour, then the appropriate 1-handles and 0-handle together with the neighbourhood $D$ of T-center construct an oriented surface. The change of the orientation of the boundary of $D$ and the colours of all three edges give us an oriented surface too, but all the other cases cause appearance of a M$\ddot{\textrm{o}}$bius strip.
\end{proof}

\begin{proof} (theorem \ref{t7})
If we fix an orientation in the unit circle of a T-diagram, then according to the described procedure a code is defined clearly. We get a symmetric code for another orientation. So the codes are either identical or symmetric for topologically equivalent flows. On the other hand, we can renew the chord diagram definitely with the help of the code, if we have a regular 2n+1-gon with its vertices in the unit circle (n is a quantity of mutually distinct numbers in the code apart from 0) with put numbers on its vertices according to the code (doing a clockwise motion along the circle). The pairs of the points with the same numbers be joined by the chords, whereas the center of the chord 1-1 be joined with 0 constructing T. The chords (the edges) having ' near their numbers be coloured with the second colour, all the others --- with the first one.

\end{proof}


\begin{proof} (theorem \ref{t8})
When constructing a reverse diagram, the cycle maps into a unit circle as it is a boundary of the neighbourhood of the source. Curves which are transversal to the chords and the lower edge of T are isotopic to the unstable manifolds of the singular points which correspond the couple of points in the T-diagram of the reverse flow. Coherence of orientations along these transversal curves defines a colour of chords as when constructing a T-diagram. According to the construction the T-graph of the reverse diagram is got from the T-graph of the initial diagram with an axial symmetry around a middle perpendicular to the lower edge (when an appropriate embedding both graphs into the plane). This mapping reverts an orientation. So, if the direction of motion on the unit circle of the T-diagram was defined with an ordered triple 0, 1a, 1b, then the direction on the unit circle of the reverse diagram which is coherent with it is defined with the appropriate triple 0, 1b, 1a. It defined a motion along the cycle and allows to define a colour of the edges of the T-graph.
\end{proof}

\section{Examples} \label{examples}

Since now we consider everywhere for the singularity ``saddle--node{}'' that the source is a node.

\subsection*{SN-flows on oriented surfaces}
There is a unique optimal flow with a singularity ``saddle--node{}'' (source) on a sphere. C-diagram is an empty circle (without chords) with a marked arc in it.

\begin{figure}[ht]
\center{\includegraphics[width=0.82\linewidth]{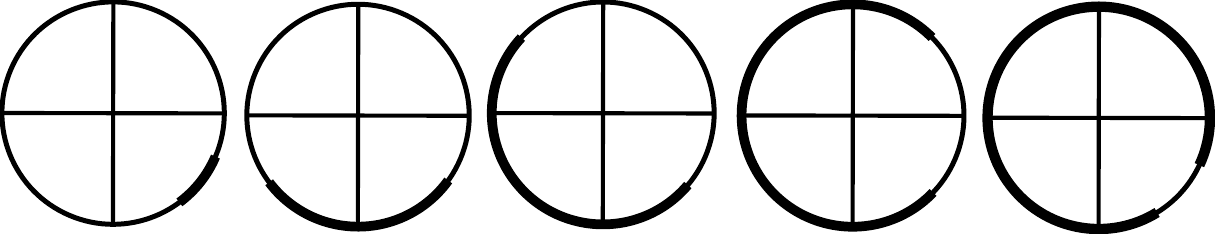}
}
\caption{SN-diagrams on the torus}
\label{bif1}
\end{figure}

There is a unique chord diagram of an optimal Morse flow on a torus. We  fix a standard orientation of a circle --- in counterclockwise direction. There is the only way of choice for the beginning of a C-arc  to within circle symmetries. Then this point divides the arc it belongs to into two parts. Hence there are five options of arcs for the end of C-arc which can lie in it. So , there are five topologically non-equivalent SN-flows (SN-node is a source) on a torus. Their C-diagrams can be seen in fig.\ref{bif1}.

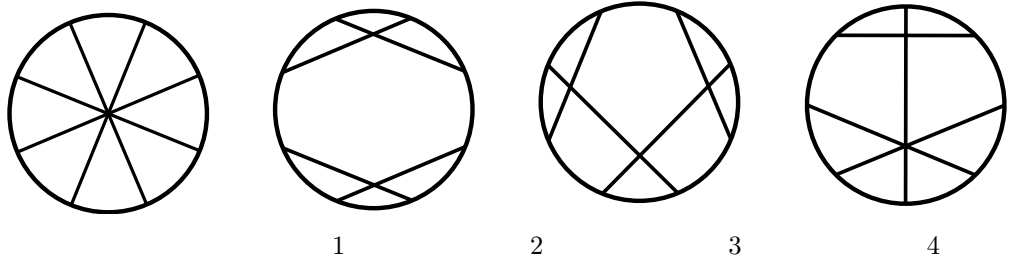
\begin{figure}[ht]
\center{\includegraphics[width=0.8\linewidth]{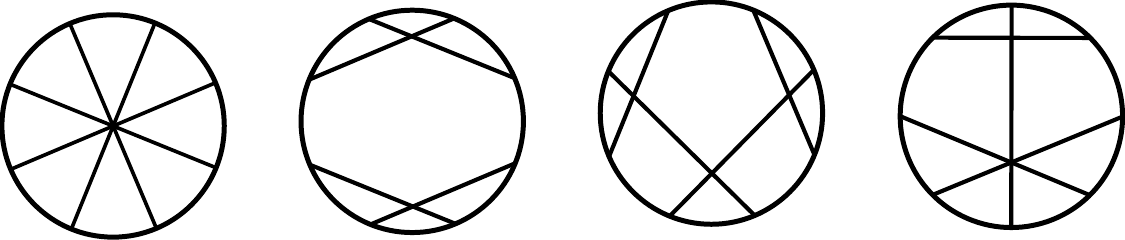}
}
\put (-255,-15) {1}
\put (-180,-15) {2}
\put (-105,-15) {3}
\put (-30,-15) {4}
\caption{Chord diagrams on the genus 2 surface}
\label{sn2}
\end{figure}

There are four chord diagrams of an optimal flow on an oriented surface of genus 2 (fig.\ref{sn2}). The beginning of a C-arc in the first diagram can be chosen in one way, whereas the end can be chosen in nine ways. So, there are nine different C-diagrams.

We can choose the first end in the second diagram in three ways: 1) right central arc; 2) right upper arc; 3) upper central arc. In the first case we have nine options for the second end, in the second case --- seven options which were not considered before, in the third case --- three options. Altogether we have nineteen opportunities for a placement of the C-arc in the second diagram.

We have five options for the first end of the C-arc in the third diagram, starting from the lower arc and doing counterclockwise circumvention. There be $9 + 8 + 6 + 4 + 2 = 29$ opportunities for the placement of this arc in general, they define twenty-nine C-diagrams.

We have four options for the first end in the fourth diagram and $9 + 7 + 5 + 3 = 24$ for the second end, so we have twenty-four C-diagrams in this case.

So, in general we have next statement: there are eighty-one SN-flows (source) on an oriented surface of genus 2.

\subsection*{SN-flows on non-oriented surfaces}
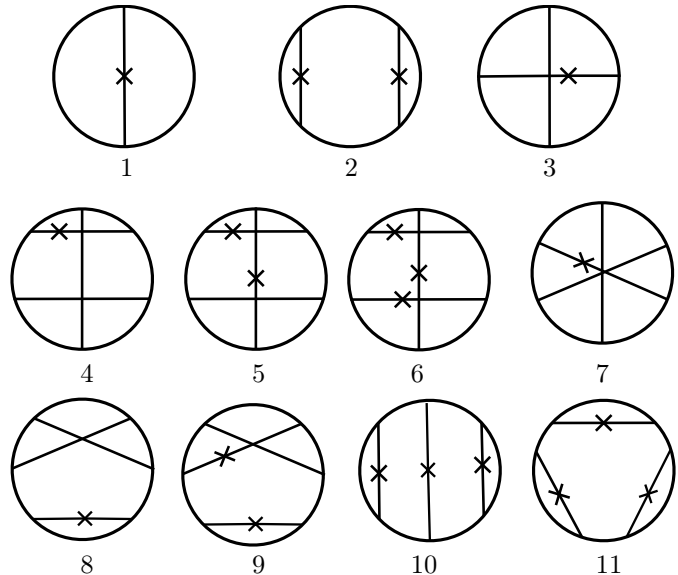
\begin{figure}[ht]
\center{\includegraphics[height=7.2cm]{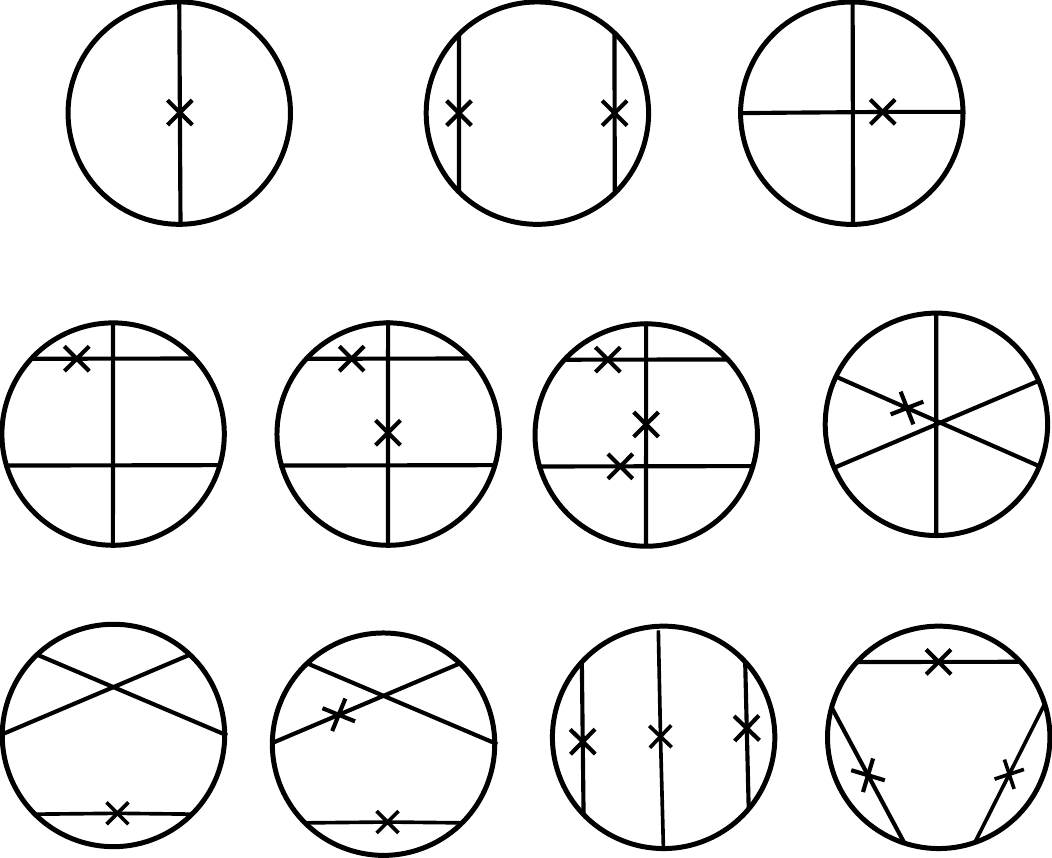}}
\put (-210,140) {1}
\put (-125,140) {2}
\put (-50,140) {3}
\put (-225,62) {4}
\put (-160,62) {5}
\put (-100,62) {6}
\put (-30,62) {7}
\put (-225,-10) {8}
\put (-160,-10) {9}
\put (-100,-10) {10}
\put (-30,-10) {11}
\caption{Non-oriented chord diagrams:
\\    1 -- genus 1,  \ \
 2,3 -- genus 2, \ \ 4--11 -- genus 3.}
\label{snn}
\end{figure}

Chord diagrams of optimal Morse flows on non-oriented surfaces of genuses 1, 2 and 3 (\cite{Kybalko2018}) are in the fig. \ref{snn}.

We have the first diagram on the projective plane (genus 1). The first end of the C-arc on it is defined clearly, whereas there are three options for the second one. So, there are three topologically unequivalent SN-flows (source) on the projective plane.

We have two diagrams (2 and 3 in the fig. \ref{snn}) on the Klein bottle (genus 2).
There are two options for the first end of the C-arc in the diagram 2, 5+3=8 options for the second one. The first end of the C-arc is defined clearly for the diagram 3, but we have five opportunities for the second one. Hence there are thirteen topologically unequivalent SN-flows (source) on the Klein bottle.

There are eight chord diagrams on a non-oriented surface of genus 3: diagrams 4--11 in the fig. \ref{snn}. Let's find options of placement of the C-diagram according to the previous scheme.  We are getting next facts:

for diagrams 4 and 5 number of options is equal to 7+5+3=15 for each one,

for diagrams 6 and 7 number of options is equal to 7+5=12 for each one,

for diagram 8 number of options is equal to 7+6+4+2=19,

for diagram 9 number of options is equal to 7+6+5+4+3+2=27,

for diagram 10 number of options is equal to 7+5=12,

for diagram 11 number of options is equal to 7+4=11.

Altogether we have one hundred twenty-three flows on a non-oriented surface of genus 3.








\subsection*{SC-flows on oriented surfaces}
Let's suppose we have a sphere and an optimal SC-flow on it. A saddle connection joins two saddles. Their Poincar$\acute{\textrm{e}}$ rotation index is equal to $-1$ for each one. Then, according to the Poincar$\acute{\textrm{e}}$--Hopf theorem, general number of all the sources and sinks is equal to four. Three cases are possible: 1) one source and three sinks, 2) two sources and two sinks, 3) three sources and one sink. Everywhere we have a T-graph with no chords.

In the first case one unit circle corresponds to one source. All the edges of the T-graph have the first colour. In the second case we have two circles and the vertices of the T-graph lie in them. Two subcases are possible: a) two side vertices lie in one circle, whereas the lower vertice does in the other one, b) side vertices lie in different circles. In the third case we have three circles and each one of them has one vertice of the T-graph in it. So, general number of optimal SC-flows on a sphere is equal to four.

There is a unique optimal SC-flow on a torus, because there is a unique T-diagram without chords which corresponds a torus --- all the edges are coloured with the second colour.

\begin{figure}[ht!]
\center{\includegraphics[width=0.7\linewidth]{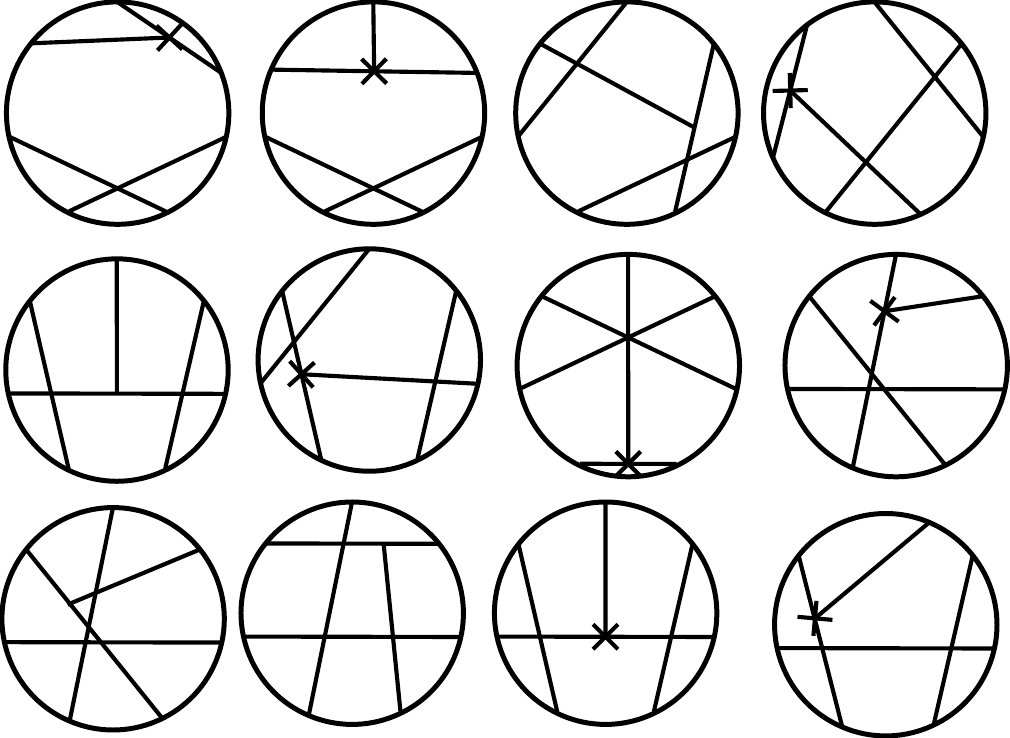}
}
{
\caption{Twelve saddle connection bifurcations on the genus 2 surface}
}
\label{sc2}
\end{figure}

The T-diagram of an optimal SC-flow on an oriented surface of genus 2 has two chords and a T-graph. All such T-diagrams which are possible are showed in the fig. \ref{sc2}. Here, a graph $T$ with a mark at the central vertex denotes a graph in which all three edges are marked.

\subsection*{SC-flows on non-oriented surfaces}
\begin{figure}[ht]
\center{\includegraphics[height=2cm]{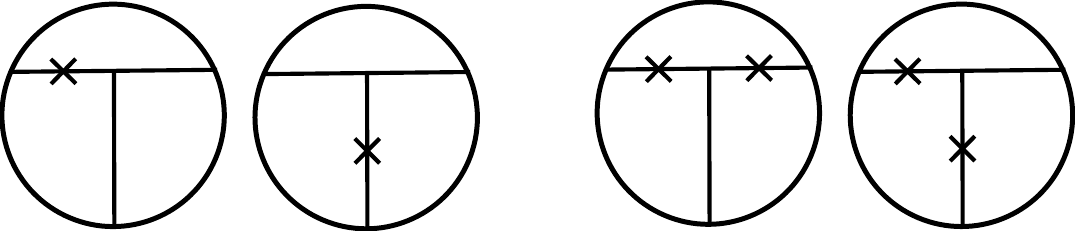}}
\put (-240,-15) {1}
\put (-175,-15) {2}
\put (-92,-15) {3}
\put (-27,-15) {4}
\caption{Non-oriented T-diagrams without chords}
\label{T2}
\end{figure}
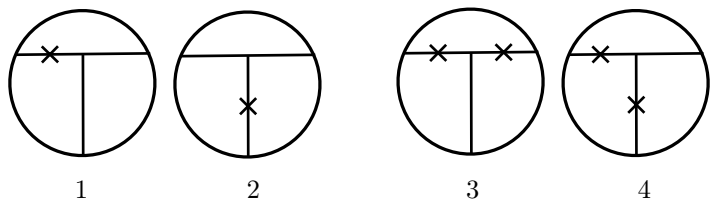

If a flow with hyperbolic singular points on the projective plane has two saddles, then it has yet three points which are sources and sinks. We are analyzing flows with one source and two sinks. The T-diagram of such flows has no chords and has two cycles. We can see such diagrams in fig. \ref{T2}.1 and 2.

Flows with two sources and one sink we can get from these flows with reversing orientations of all trajectories. So we have four topologically unequivalent optimal SC-flows on the projective plane.

An optimal SC-flow on the Klein bottle has one source, one sink and two saddles. So its T-diagram has no chords and has one cycle. Two such diagrams which are possible we can see in \ref{T2}.3 and 4. Hence there are two topologically unequivalent SC-flows on the Klein bottle.

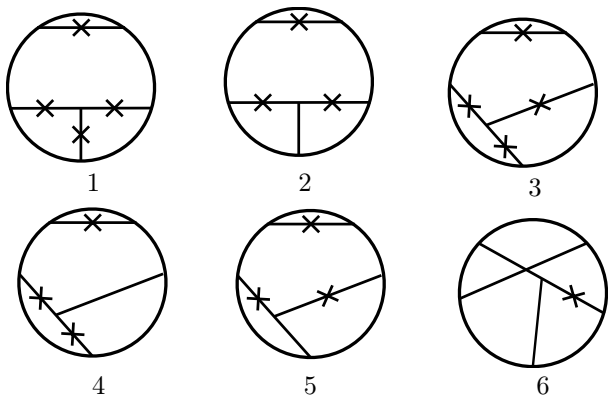
\begin{figure}[ht]
\center{\includegraphics[height=4.8cm]{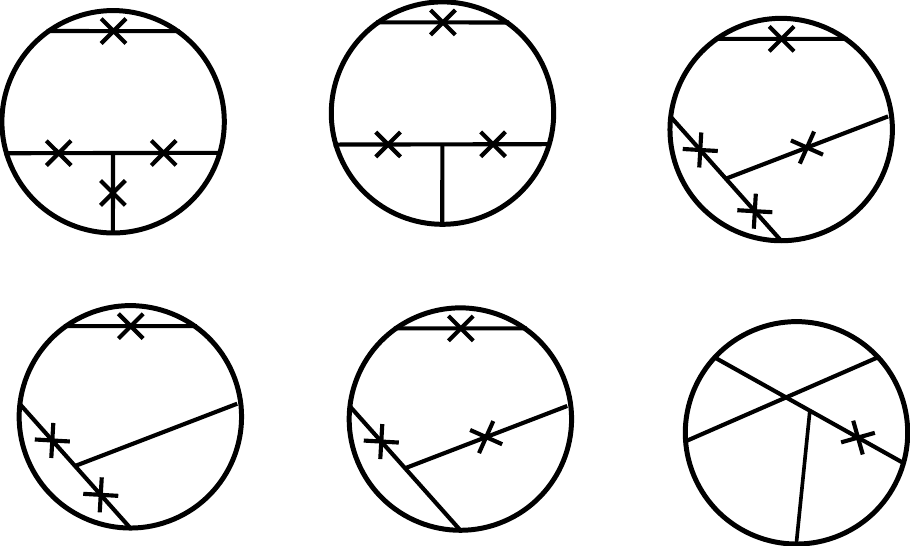}}
\put (-197,67) {1}
\put (-117,67) {2}
\put (-30,65) {3}
\put (-195,-10) {4}
\put (-115,-10) {5}
\put (-27,-10) {6}
\caption{Non-oriented self-reverse Т-diagrams of type 3}
\label{T3a}
\end{figure}

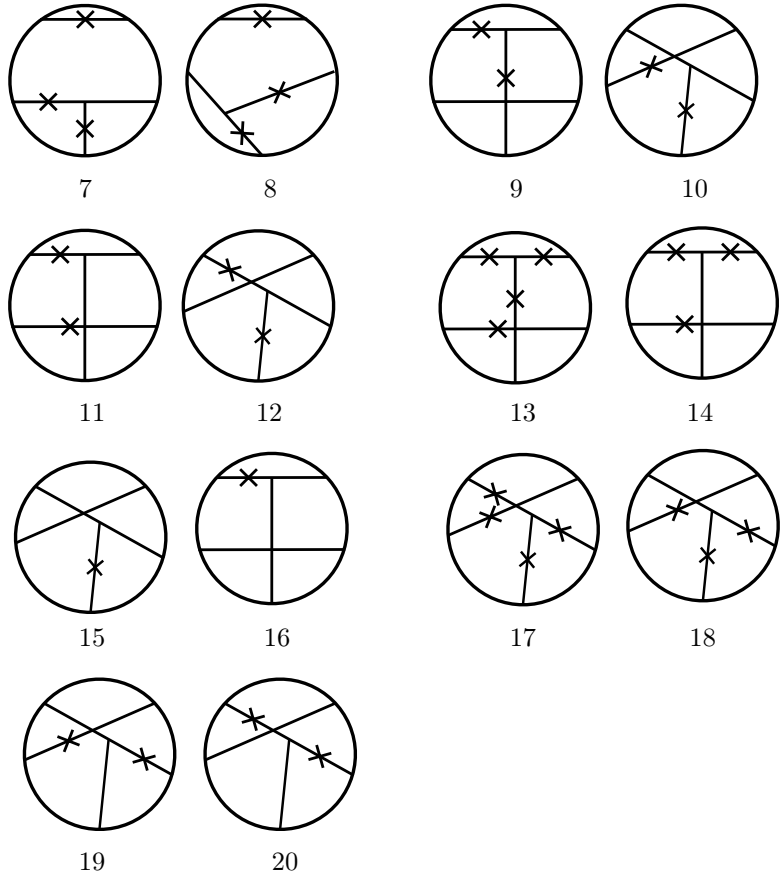
\begin{figure}[ht]
\center{\includegraphics[height=11cm]{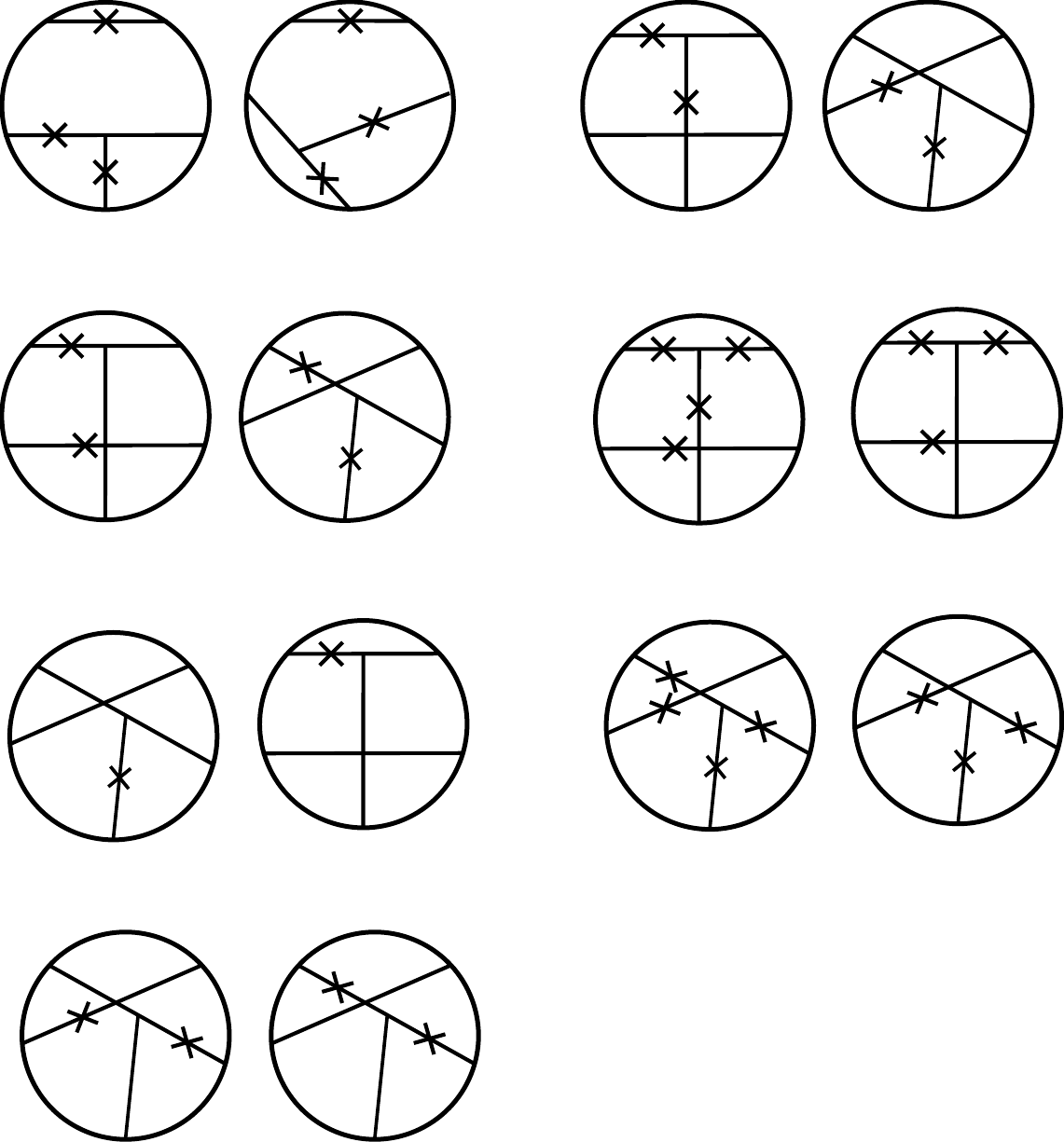}}
\put (-265,240) {7}
\put (-195,240) {8}
\put (-102,240) {9}
\put (-37,240) {10}
\put (-265,155) {11}
\put (-198,155) {12}
\put (-102,155) {13}
\put (-35,155) {14}
\put (-265,70) {15}
\put (-195,70) {16}
\put (-102,70) {17}
\put (-34,70) {18}
\put (-265,-15) {19}
\put (-192,-15) {20}

\caption{Pairs of mutually reverse non-oriented Т-diagrams of type 3}
\label{T3b}
\end{figure}

The T-diagrams of optimal SC-flows on a non-oriented surface of genus 3 can be seen in fig. \ref{T3a} and \ref{T3b}.

First we coloured the edges of the T-graph with two colours for finding these diagrams. The case when all three edges have the first colour is impossible because one can't join three constructed cycles with one chords. If only one edge is coloured with the second colour, then the chord has to join two constructed cycles which are an arc two edges of the first colour are attached to and one of the other two arcs. Colouring the chord can be arbitrary. So we have four diagrams for the side edge of the second colour: 6, 11, 16, 19. There are two diagrams for the lower edge --- 10 and 15.


If two or three edges are coloured with the second colour, then the chord can join two arbitrary points in the circle, but only one out of two colourings gives one cycle in the diagram in this case.

If the side edge has the first colour, then six options are possible for the chord: three options with a chord on one arc --- 5, 7, 8, and three options with chords on different arcs -- 9, 12, 18.

If the lower edge has the first colour, but the side ones have the second colour, then we have two options when the ends of the chord belong to the same arc --- 2 and 4, and two options when they belong to different ones --- 14 and 20.

If all the edges have the second colour, then the chord has it too. For the chord on one arc we have two options --- 1 and 3, otherwise we also have two options --- 13 and 17.

So, there are all possible twenty non-oriented diagrams with one edge in fig. \ref{T3a} and \ref{T3b}.

We can use the method described in the previous chapter to verify which diagrams here are mutually reverse.

Hence there are twenty topologically unequivalent SC-flows on a non-oriented surface of genus 3, and there are six ones among them reverse to which are topologically equivalent to themselves and yet seven pairs of mutually reverse flows.

\section{Conclusion}

In this paper  we construct complete topological invariants of codimension one gradient flows on closed surfaces. They were used for description of all the possible structures of flows on oriented surfaces of genus 0, 1, 2 and non-oriented surfaces of genus 1, 2, 3. Also we described an approach of finding the invariants of reverse flows according to the invariant of the initial flow.

We hope that results we have got here will be generalized on the surfaces with boundary and on arbitrary codimension one flows.



\end{document}